\documentclass[oneside,notitlepage,12pt]{article}

\usepackage{amssymb}
\usepackage[leqno]{amsmath}
\usepackage{amsfonts}
\usepackage{amsopn}
\usepackage{amstext}
\usepackage{amsthm}

\usepackage{verbatim}
\usepackage[colorlinks, backref]{hyperref}
\usepackage{makeidx}

\usepackage{calrsfs}
\usepackage{fourier-orns}
\usepackage{clock} 
\usepackage[alpine, weather]{ifsym}

\providecommand{\cal}{\mathcal}

\newenvironment{pf}{\begin{proof}}{\end{proof}}

\newcommand{\Aaa}{{\cal{A}}}
\newcommand{\Bee}{{\cal{B}}}
\newcommand{\Cee}{{\cal{C}}}
\newcommand{\Dee}{{\cal{D}}}

\newcommand{\Tau}{{\cal{T}}}
\newcommand{\Yu}{{\cal{U}}}

\newcommand{\Wu}{{\cal{W}}}

\newcommand{\al}{\alpha}

\renewcommand{\phi}{\varphi}
\renewcommand{\rho}{\varrho}

\newcommand{\ntr}{{n\in\omega}}

\newcommand{\loe}{\leq}
\newcommand{\goe}{\geq}

\newcommand{\subs}{\subseteq}
\newcommand{\sups}{\supseteq}
\newcommand{\nnempty}{\ne\emptyset}

\renewcommand{\iff}{\Longleftrightarrow}

\newcommand{\id}[1]{{\operatorname{id}_{#1}}} 

\newcommand{\defi}{\,:=\,}

\newtheorem{tw}{Theorem}[section]
\newtheorem{wn}[tw]{Corollary}
\newtheorem{lm}[tw]{Lemma}
\newtheorem{prop}[tw]{Proposition}

\theoremstyle{definition}

\newtheorem{ex}[tw]{Example}
\newtheorem{exs}[tw]{Examples}

\theoremstyle{remark}

\newcommand{\setof}[2]{\{#1\colon #2\}}
\newcommand{\bigsetof}[2]{\Bigl\{#1\colon #2\Bigr\}}

\newcommand{\sett}[2]{\{#1\}_{#2}}
\newcommand{\sn}[1]{\{#1\}} 
\newcommand{\dn}[2]{\{#1,#2\}} 
\newcommand{\map}[3]{#1\colon #2 \to #3} 

\newcommand{\bF}{{\mathbb F}}

\newcommand{\calN}{\mathcal N}

\newcommand{\ci}{\operatorname{ci}}
\newcommand{\cir}{\operatorname{cir}}

\newcommand{\cB}{\mathcal B}
\newcommand{\bfind}[1]{\index{#1}{\em #1}}

\title{Chain intersection closures}
\author{
{\sc Wies{\l}aw Kubi\'s}\footnote{The first author was supported by GA\v CR grant No. GA17-27844S (Czech Science Foundation).} \\
{\small Mathematical Institute, Czech Academy of Sciences, Czech Republic}\\
{\small Cardinal Stefan Wyszy\'nski University, Warsaw, Poland}
\and
{\sc Franz-Viktor Kuhlmann}\footnote{The second author was partially supported by a Polish Opus grant 2017/25/B/ST1/01815.} \\
{\small Institute of Mathematics, University of Szczecin, Poland}
}

\date{October 7, 2018}

\begin{document}

\maketitle

\begin{abstract}
	We study spherical completeness of ball spaces and its stability under expansions. We introduce the notion of an ultra-diameter, mimicking diameters in ultrametric spaces. We prove some positive results on preservation of spherical completeness involving ultra-diameters with values in narrow partially ordered sets. Finally, we show that in general, chain intersection closures of ultrametric spaces with partially ordered value sets do not preserve spherical completeness.


\end{abstract}


\section{Introduction}

In \cite{KK1,KK1a,KK2,KK3,KK4,KK5}, ball spaces are studied in order to provide a general framework for fixed
point theorems that in
some way or the other work with contractive functions. A {\em ball space} $(X,\cB)$ is a nonempty set $X$ together
with any nonempty collection of nonempty subsets of $X$. The completeness property necessary for the proof of
fixed point theorems is then encoded as follows. A \bfind{chain of balls} (also called a \emph{nest})
in $(X,\cB)$ is a nonempty subset of $\cB$ which is linearly ordered by inclusion. A ball space $(X,\cB)$ is
called \bfind{spherically complete} if every chain of balls has a nonempty intersection. Further, we
say that a ball space $(X, \cB)$ is \bfind{chain intersection closed} if the intersection of every chain in
$\cB$ is either empty or a member of $\cB$. We define $\ci(\cB)$ to be the family of all nonempty sets of the form $\bigcap \Cee$, where $\Cee \subs \Bee$ is a chain (recall that, by default, chains of sets are supposed to be nonempty). More formally,
$$\ci(\cB) = \bigsetof{ \bigcap \Cee }{\emptyset \ne \Cee \subs \cB, \text{ $\Cee$ is a chain} } \setminus \sn \emptyset.$$
Hence a ball space $(X, \cB)$ is chain intersection closed if and only if $\ci(\cB) = \cB$. In the present paper, we study the process of obtaining a chain intersection closed ball space from a given ball space and the question under which conditions the spherical completeness of $(X,\cB)$ implies the spherical completeness of $(X,\ci(\cB))$.

Main inspiration for these definitions and questions is taken from the theory of ultrametric spaces and their ultrametric balls.
An \bfind{ultrametric} $d$ on a set $X$ is a function from $X\times X$ to a partially ordered set
$\Gamma$ with smallest element $\bot$, such that for all $x,y,z\in X$ and all $\gamma\in\Gamma$,
\begin{enumerate}
	\item[(U1)] $d(x,y)=\bot$ if and only if $x=y$,
	\item[(U2)] if $d(x,y)\leq\gamma$ and $d(y,z)\leq\gamma$, then
	$d(x,z)\leq\gamma$,
	\item[(U3)] $d(x,y)=d(y,x)$ \ \ \ (symmetry).
\end{enumerate}
Condition (U2) is the ultrametric triangle law; if $\Gamma$ is linearly ordered, it
can be replaced by
\begin{enumerate}
	\item[(UT)] $d(x,z)\leq\max\{d(x,y),d(y,z)\}$.
\end{enumerate}
A \bfind{closed ultrametric ball} is a set $B_\alpha(x) \defi \setof{y\in X}{d(x,y)\leq \alpha}$, where $x\in X$ and $\alpha\in\Gamma$. The problem with general ultrametric spaces is that closed balls $B_\alpha(x)$ are not necessarily precise, that is, there may not be any $y\in X$ such that $d(x,y)=\alpha$. Therefore, we prefer to work only with \bfind{precise} ultrametric balls, which we can write in the form
$$
B(x,y)\defi \setof{z\in X}{d(x,z)\leq d(x,y)},
$$
where $x,y\in X$. We obtain the \bfind{ultrametric ball space} $(X,\cB_d)$ from $(X,d)$
by taking $\cB_d$ to be the set of all such balls $B(x,y)$.
Specifically, $\cB_d \defi \setof{B(x,y)}{x,y \in X}$.

More generally, an ultrametric ball is a set $B_S (x) \defi \setof{y\in X}{d(x,y)\in S}$, where $x\in X$ and $S$ is
an initial segment of $\Gamma$. We call $X$ together with the collection of all ultrametric balls
the \bfind{full ultrametric ball space} of $(X,d)$. Any nonempty intersection of a chain $(B(x_i,y_i))_{i\in I}$
is such an ultrametric ball $B_S (x)$, where $S$ is the initial segment $\setof{\gamma\in\Gamma}{\gamma\leq
d(x_i,y_i)\mbox{ for all }i\in I}$. If $\Gamma$ is linearly ordered, then $B_S (x)$ is the intersection of the
chain $\setof{B(x,y)}{y\in X\mbox{ such that }d(x,y)\notin S}$; hence in this case, the full ultrametric ball space is
just $(X,\ci(\cB_d))$.

\begin{tw}                                  \label{MT1}
Let $(X,\cB_d)$ be the ball space of an ultrametric space $(X,d)$ with linearly ordered value set. Then the following
assertions hold:
\par\smallskip\noindent
1) \ The ball space $(X,\ci(\cB_d))$ is chain intersection closed.
\par\smallskip\noindent
2) \ If $(X,\cB_d)$ is spherically complete, then so is $(X,\ci(\cB_d))$.
\end{tw}
\noindent
We will deduce this theorem from the more general Theorem~\ref{ThmDwaJedna} in the next section.

\par\smallskip
Assume that $\setof{B_i}{i\in I}$ is any collection of balls in $\ci(\cB_d)$ such that $B_i\cap B_j\ne\emptyset$
for all $i,j\in I$. Then from the ultrametric triangle law and the assumption that the value set is linearly
ordered it follows that $\setof{B_i}{i\in I}$ is in fact a chain.
Hence it follows from part 1) of our theorem that $\ci(\cB_d)$ is closed under nonempty intersections of
{\em arbitrary} collections of balls.

\par\medskip
The structure of ultrametric spaces with partially ordered value sets is in general much more complex than in the
case of linearly ordered value sets. But we can at least prove the following.
Recall that a partially ordered set (``poset'') is \bfind{narrow} if it contains no infinite sets of pairwise
incomparable elements.

\begin{tw}                                     \label{MT2}
Let $(X,\cB_d)$ be the ball space of an ultrametric space with countable narrow value set. Then the assertions of
Theorem~\ref{MT1} hold.
\end{tw}
\noindent
This theorem will be proven at the end of Section~\ref{Sectud}, where we study intersection closures for the more
general class of ball spaces that admit functions which we call ``ultra-diameters''. These are functions that
associate to every ball a value in a poset having a special property related to ultrametrics. For instance, for a
ball space consisting of precise balls in an ultrametric space $(X,d)$, the function $B(x,y)\mapsto d(x,y)$ is an
ultra-diameter.

\par\smallskip
Take two ball spaces $(X,\cB)$ and $(X,\cB')$ on the same set $X$. We call $(X,\cB')$ an \bfind{expansion}
of $(X,\cB)$ if $\cB\subseteq \cB'$. In general, we cannot expect the existence of chain
intersection closed expansions which preserve spherical completeness. Example~\ref{Ex1} in
Section~\ref{Sectnr} shows that the condition ``narrow'' cannot be dropped in Theorem~\ref{MT2}:

\begin{tw}                                     \label{MT3}
There exists a countable spherically complete ultrametric space with partially ordered value set whose
ultrametric ball space does not admit any expansion that is chain intersection closed and spherically complete.
\end{tw}

However, we do not know whether the countability assumption can be dropped in Theorem~\ref{MT2}. Example~\ref{Ex2} 
in Section~\ref{Sectnr} presents an uncountable narrow spherically complete ball space which does not admit any 
expansion that is chain intersection closed and spherically complete.
But in contrast to Example~\ref{Ex1}, it cannot be transformed into an ultrametric space.

\section{Chain intersection closure}
Let $\cB$ be a nonempty family of nonempty sets. Using transfinite recursion, we define $\ci_\al(\cB)$ for each
ordinal $\al$, as follows.
$$\ci_0(\cB) = \cB, \qquad \ci_\al(\cB) = \ci\left( \bigcup_{\xi < \al} \ci_\xi(\cB) \right) \text{ for } \al > 0.$$
Finally, we define the \bfind{chain intersection rank} of $\cB$, denoted by $\cir(\cB)$, to be the smallest ordinal $\al$ such that $\ci_{\al+1}(\cB) = \ci_\al(\cB)$.
Thus, $\cir(\cB) = 0$ if and only if $\cB$ is chain intersection closed, while $\cir(\cB) \loe 1$ means that in order to make $\cB$ chain intersection closed, it suffices to extend it by adding all nonempty intersections of chains.
In general, we call $(X, \ci_\al(\cB))$, with $\al = \cir(\cB)$, the \bfind{chain intersection closure} of $(X, \cB)$. It could also be described as a ball space $(X, \cB')$, where $\cB' \sups \cB$ is minimal such that $\cB' \cup \sn \emptyset$ is stable under intersections of chains.

A ball space $(X,\cB)$ will be called \bfind{tree-like} if for every $B_1, B_2 \in \cB$ the following implication
holds.
\begin{equation}
	B_1 \cap B_2 \nnempty \implies B_1 \subs B_2 \text{ or } B_2 \subs B_1.
	\tag{I}\label{eqaI}
\end{equation}
Note that if $(X, \cB)$ is tree-like then the poset $(\cB, \sups)$ is a \bfind{generalized tree} in the sense that for every $B \in \cB$ the set $\setof{C \in \cB}{C \sups B}$ is linearly ordered. The converse may be false, as the relation $B_1 \cap B_2 \nnempty$ does not imply the existence of any $B \in \cB$ with $B \subs B_1 \cap B_2$.
Perhaps the simplest counterexample is $\cB = \{ \{0,1\}, \{1,2\} \}$, where $(\cB, \sups)$ is a tree consisting of two incomparable elements, however (\ref{eqaI}) is violated. But the following is apparent form our discussion:
\begin{lm}
Let $(X, \cB)$ be a ball space such that $\cB$ is closed under finite intersections. If $(\cB, \sups)$ is a
generalized tree, then $(X, \cB)$ is a tree-like ball space.
\end{lm}

This is our main theorem on tree-like ball spaces.
\begin{tw}\label{ThmDwaJedna}
	Let $(X, \cB)$ be a tree-like ball space. Then
	\begin{enumerate}
		\item[{\rm(1)}] $\cir(\cB) \loe 1$,
		\item[{\rm(2)}] $(X,\ci(\cB))$ is tree-like,
		\item[{\rm(3)}] $(X,\ci(\cB))$ is spherically complete whenever $(X,\cB)$ is.
	\end{enumerate}
	
\end{tw}

\begin{pf}
	Let $\Dee$ be a chain in $\ci(\cB)$.
	For each $D \in \Dee$, choose a chain $\Cee_D \subs \cB$ such that $D = \bigcap \Cee_D$.
	Let $\Cee = \bigcup_{D \in \Dee}\Cee_D$.
	Then $\bigcap \Cee = \bigcap \Dee$.
	We claim that $\Cee \subs \cB$ is a chain. Indeed, fix $C_1, C_2 \in \Cee$ and consider $D_1, D_2 \in \Dee$ such that $C_i \in \Cee_{D_i}$ for $i=1,2$.
	Since $\Dee$ is a chain, we may assume that $D_1 \subs D_2$ (the other possibility is the same).
	Now $C_1 \cap C_2 \sups D_1$; since $D_1$ is nonempty, we may apply (\ref{eqaI}) to obtain that either $C_1 \subs C_2$ or $C_2 \subs C_1$.
	This proves (1) and (3).
	
	In order to show (2), fix $D_i = \bigcap \Cee_i$, $i=1,2$, where $\Cee_1$, $\Cee_2$ are chains in $\Bee$ and suppose $D_1 \cap D_2 \nnempty$.
	Then $C_1 \cap C_2 \nnempty$ for every $C_1 \in \Cee_1$, $C_2 \in \Cee_2$, therefore by (\ref{eqaI}) the family $\Cee_1 \cup \Cee_2$ is a chain.
	Now, if $\Cee_1$ is co-initial in $(\Cee_1 \cup \Cee_2, \subs)$ then $D_1 \subs D_2$. Otherwise, $\Cee_2$ must be co-initial in $(\Cee_1 \cup \Cee_2, \subs)$, yielding $D_2 \subs D_1$. This shows (2).
\end{pf}

Every ultrametric space with linearly ordered value set is tree-like, since in this case property (\ref{eqaI})
follows from the ultrametric triangle law. Hence Theorem~\ref{MT1} is a special case of Theorem~\ref{ThmDwaJedna}.

\par\medskip
The following example shows that the chain intersection rank of a ball space can be arbitrarily large.
\begin{ex}\label{ExNoWERTwegho}
	Let $\Aaa$ be an arbitrary family of nonempty sets.
	We claim that there is a family of sets $\Bee$ which is a generalized tree both with $\subs$ and with $\sups$
    and such that $\ci(\Bee) = \Bee \cup \Aaa$.
	This will also show that the chain intersection rank can have arbitrarily large values.
	
	For each $A \in \Aaa$ choose a countable set $E(A) = \sett{e_{A,n}}{\ntr}$ so that $E(A) \cap E(B) = \emptyset$ whenever $A \ne B$.
	Define $E_n(A) = A \cup \setof{e_{A,i}}{i \goe n}$.
	Let $$\Bee = \setof{E_n(A)}{\ntr, \; A \in \Aaa}.$$
	Clearly, $\sett{E_n(A)}{\ntr}$ is a chain satisfying $\bigcap_{\ntr}E_n(A) = A$ for every $A \in \Aaa$, therefore $\Aaa \subs \ci(\Bee)$.
	On the other hand, $E_k(A)$ and $E_\ell(B)$ are comparable with respect to inclusion if and only if $A = B$, therefore $\ci(\Bee) = \Bee \cup \Aaa$.
	Using transfinite induction, we can repeat this construction as long as we wish, obtaining families of sets with arbitrarily large chain intersection ranks.
\end{ex}

The result above motivates the following definition. Namely, we say that a ball space $(X,\cB)$ is \bfind{chain intersection stable} if for every nonempty family $\bF$ consisting of chains in $\cB$ such that $\setof{\bigcap \Cee}{\Cee \in \bF}$ is a chain, there exists a chain $\Yu \subs \cB$ satisfying
$$\bigcap \Yu = \bigcap_{\Cee \in \bF} \bigcap \Cee.$$
Clearly, if $(X, \cB)$ is chain intersection stable then $\cir(\cB) \loe 1$ and spherical completeness is preserved
when passing from $\cB$ to $\ci(\cB)$. In order to prove that a concrete ball space $(X,\cB)$ is chain intersection stable, some diagonalization argument needs to be invoked. In the proof of Theorem~\ref{ThmDwaJedna}, simply the union of all chains is a chain. So we have:

\begin{wn}                                   \label{tlcis}
Every tree-like ball space is chain intersection stable. In particular, every ultrametric space with linearly ordered value set is chain intersection stable.
\end{wn}

On the other hand, a ball space with chain intersection rank $\leq 1$ is not necessarily chain intersection stable.
Indeed, if in Example~\ref{ExNoWERTwegho} there is a chain in $\Aaa$ with nonempty intersection not contained in
$\Aaa$, then the intersection over this chain is not equal to the intersection over any chain in $\Bee$.

We shall see in the proof of Theorem~\ref{ThmHevyNarow} that sometimes a quite nontrivial diagonalization is
needed to prove that a given ball space is chain intersection stable.

\section{Ultra-diameters}                      \label{Sectud}
Recall that a function $\map{f}{(P,\loe_P)}{(Q,\loe_Q)}$ is \bfind{order preserving} or \bfind{increasing} if $x_0 \loe_P x_1$ implies $f(x_0) \loe_Q f(x_1)$, while it is called \bfind{strictly order preserving} or \bfind{strictly increasing} if additionally $f(x_0) \ne f(x_1)$ whenever $x_0 \ne x_1$.
Below we define the concept of an ultra-diameter, from a family of sets to a fixed poset. Note that every family of sets is a poset with inclusion or reversed inclusion. On the other hand, every poset is isomorphic to a family of sets with inclusion (as well as reversed inclusion).

An \bfind{ultra-diameter} on a family of nonempty sets $\cB$ is a function $\map \delta \cB \Gamma$, where $\Gamma$ is a poset, satisfying
\begin{enumerate}
	\item[(D1)] $\delta$ is increasing (i.e., $\delta(B_0) \loe \delta(B_1)$ whenever $B_0 \subs B_1$),
	\item[(D2)] if $B_0, B_1 \in \cB$, $B_0 \cap B_1 \nnempty$, and $\delta(B_0) \loe \delta(B_1)$ then $B_0 \subs B_1$.
\end{enumerate}

Of course, the identity function $\map{\id \cB}{\cB}{\cB}$ is an ultra-diameter. In general, the idea is to find $\Gamma$ as simple as possible, so that there is still an ultra-diameter from $\cB$ to $\Gamma$.
The adjective ``ultra" is motivated by (generalized) ultrametric spaces in which the precise balls have a natural ultra-diameter.

The proof of the following lemma is straightforward, by applying (U2) and (U3).

\begin{lm}\label{LmXxxWERTUGuwrg}
	Let $(X,d)$ be an ultrametric space with value poset $\Gamma$. Let $x,y \in X$ be such that $\gamma = d(x,y)$.
	Then $B(x,y) = B(x', y')$ for every $x', y' \in B(x,y)$ satisfying $d(x', y') = \gamma$.
\end{lm}

By Lemma~\ref{LmXxxWERTUGuwrg}, the function $\map{\delta}{\cB_d}{\Gamma}$ given by the formula
$$\delta(B(x,y)) = d(x,y)$$
is well defined.
It obviously satisfies (D1).
If $B(x,y) \cap B(u,w) \nnempty$ and, say, $d(x,y) \loe d(u,w)$ then by the ultrametric triangle law (U2) we get
$B(x,y) \subs B(u,w)$, therefore $\delta$ satisfies (D2), showing that it is an ultra-diameter on $\cB_d$.

Not every ultra-diameter comes from an ultrametric, simply because an ultrametric is defined on points, hence on a countable set it attains countably many values only. On the other hand, one can have an uncountable family $\cB$ of subsets of a fixed countable set, with an ultra-diameter $\map{\delta}{\cB}{\Gamma}$ where $\Gamma$ is uncountable and $\delta$ attains all possible values (in fact, $\Gamma$ could be $\cB$ and $\delta$ could be the identity).

\par\medskip
Recall that a poset is \bfind{narrow} if it contains no infinite sets of pairwise incomparable elements. Besides 
linearly ordered sets, finite products of ordinals endowed with the coordinate-wise ordering provide examples of 
narrow posets. For our next theorem, we will need two lemmas that reflect important and well known properties of 
narrow posets. Recall that a family $\Dee$ of sets is \bfind{linked} if $B \cap B' \nnempty$ for any $B, B' \in 
\Dee$.

\begin{lm}\label{LmDusek}
Let $\map \delta \Bee \Gamma$ be an ultra-diameter such that $\Gamma$ is narrow. Then every infinite linked family 
$\Dee \subs \Bee$ contains a chain $\Cee$ such that $|\Cee| = |\Dee|$.	Let $(P,\loe)$ be a narrow poset, 
$A \subs P$ infinite. Then there exists a chain $C \subs A$ such that $|C| = |A|$.
\end{lm}

\begin{pf}
	Given a pair $\dn B{B'}$ in $\Dee$ color it \emph{green} if $B \subs B'$ or $B' \subs B$. Otherwise, color it \emph{red}.
	By the Erd\H os-Dushnik-Miller Theorem, there is $\Cee \subs \Dee$ of the same cardinality as $\Dee$ and such that all pairs in $\Cee$ have the same color.
	Suppose that this color is red.
	Then $\delta(B)$, $\delta(B')$ are incomparable in $\Gamma$, because of the definition of an ultra-diameter and the fact that $B \cap B' \nnempty$. But $\Gamma$ is narrow, so the set $\setof{\delta(B)}{B \in \Cee}$ is finite, a contradiction.
	Thus the color of every pair of balls in $\Cee$ is green, that is, $\Cee$ is a chain.
\end{pf}

Recall that a set $A$ in a poset $(P,\loe)$ is \bfind{directed} if for every $a_0,a_1 \in A$ there is $b \in A$ with $a_0 \loe b$ and $a_1 \loe b$.
It turns out that every directed narrow poset contains a cofinal subset isomorphic to a finite product of regular cardinals, see~\cite{Fraisse}.
The following fact can also be found in~\cite{Fraisse}. We give a proof for the sake of completeness.

\begin{lm}\label{LmFinDirek}
	Every narrow poset is a finite union of directed subsets.
\end{lm}

\begin{pf}
	Let $(P,\loe)$ be a narrow poset. Passing to a suitable cofinal subset, we may assume it is well-founded.
	Now the narrowness implies that the set of all initial segments of $P$ is well-founded. Thus, supposing $P$ is not a finite union of directed subsets, we may choose a minimal (with respect to inclusion) initial segment $I \subs P$ with the same property. In particular, $I$ is not directed, so there are $a,b \in I$ such that no $x \in I$ satisfies $a \loe x$, $b \loe x$.
	Define $I_a = \setof{x \in I}{b \not \loe x}$ and $I_b = \setof{x \in I}{a \not \loe x}$.
	Then $I_a$, $I_b$ are proper initial segments of $I$ and $I = I_a \cup I_b$.
	By minimality, $I_a$ and $I_b$ are finite unions of directed subsets, therefore so is $I$, a contradiction.
\end{pf}

\par\medskip
We will now generalize Corollary~\ref{tlcis} to a larger class of ball spaces.
\begin{tw}                                          \label{ThmHevyNarow}
Let $(X,\cB)$ be a ball space such that $\cB$ admits an ultra-diameter with values in a countable narrow poset.
Then $(X,\cB)$ is chain intersection stable.
\end{tw}

\begin{pf}
Let $\map \delta \cB \Gamma$ be an ultra-diameter such that $\Gamma$ is a countable narrow poset. Note that each chain in $\cB$ is countable. Indeed, by (D1) and (D2), if $B_0 \subs B_1$ and $\delta(B_0) = \delta(B_1)$ then $B_0 = B_1$, therefore $\delta$ is one-to-one on each chain in $\cB$.
	
	Let $\sett{\Cee_\al}{\al < \kappa}$ be a family of chains in $\cB$ such that, setting $C_\al \defi \bigcap \Cee_\al$, it holds that $C_\beta \subsetneq C_\al$ whenever $\al < \beta < \kappa$, where $\kappa$ is a fixed infinite regular cardinal.	
	We need to show that $C_\infty \defi \bigcap_{\al < \kappa} C_\al = \bigcap_{\al < \kappa} \bigcap \Cee_\al$ is the intersection of some chain in $\cB$.
	If $C_\al = \emptyset$ for some $\al<\kappa$ then $C_\infty = \emptyset$ and there is nothing to prove, so let us assume $C_\al \nnempty$ for every $\al < \kappa$.
	
	Let $\Cee_\al = \sett{B_{\al,n}}{\ntr}$, where $B_{\al,n} \sups B_{\al,m}$ whenever $n < m$.
Take any $B_{\al,n}\in \Cee_\al$ and $B_{\beta,m}\in \Cee_\beta$ with $\al\leq\beta<\kappa$. Then $\emptyset\ne
C_\beta\subseteq B_{\al,n}\cap B_{\beta,m}$. This proves that the set $\bigcup_{\al<\kappa}\Cee_\al$ is linked.
If $\bigcup_{\al<\kappa}\Cee_\al$ is finite, the assertion of our theorem is trivial, thus we may assume it is infinite. Hence by Lemma~\ref{LmDusek}, $\bigcup_{\al<\kappa}\Cee_\al$ contains a chain of cardinality  
$|\bigcup_{\al<\kappa}\Cee_\al|$. By what we have shown in the beginning, this cardinality must be countable. 

Since all $C_\al$, $\al < \kappa$, are distinct, for every $\al$ there must exist some $B_{\al,n}$ that does not 
appear in $\Cee_\gamma$ for any $\gamma<\al$. This means that $\bigcup_{\al<\kappa}\Cee_\al$ must contain at least 
$\kappa$ many balls. We conclude that $\kappa = \omega$.

	Note that $\Wu := \bigcup_{\al < \omega}\Cee_\al$ is narrow, as a poset endowed with inclusion. Indeed, suppose $\sett{W_n}{\ntr} \subs \Wu$. Then there are $k < \ell < \omega$ such that $\delta(W_k) \loe \delta(W_\ell)$ or $\delta(W_\ell) \loe \delta(W_k)$. On the other hand, $W_k \cap W_\ell \sups C_m \nnempty$ for a big enough $m$, therefore by (D2) we get $W_k \subs W_\ell$ or $W_\ell \subs W_k$.

	By Lemma~\ref{LmFinDirek}, $\bigcup_{\al < \omega}\Cee_\al = \Wu_0 \cup \dots \cup \Wu_{k-1}$, where each $\Wu_i$ is $\sups$-directed.
	Fix $\al < \omega$. For each $\ntr$ there is $i < k$ such that $B_{\al,n} \in \Wu_i$, so there is $j_\al < k$ such that $B_{\al,n} \in \Wu_{j_\al}$ for infinitely many $n$, say, $n_{\al,k}$, $k\in\omega$. It follows that
there is some $j<k$ and an infinite set $M \subs \omega$ such that $j_\al = j$ for every $\al \in M$.
	Now
$$\Cee_M \>:=\> \bigcup_{\al \in M} \setof{B_{\al,n_{\al,k}}}{k<\omega}$$
is contained in $\Wu_j$, therefore it is $\sups$-directed.
We have that $\bigcap \Cee_M = \bigcap \Cee = C_\infty$.
	Finally, since $\Cee_M$ is $\sups$-directed, it contains a chain $\Dee\subseteq\cB$ that is cofinal in $(\Cee_M, \sups)$, meaning that for each $B \in \Cee_M$ there is $D \in \Dee$ satisfying $B \sups D$.
	Hence $\bigcap \Dee = \bigcap \Cee_M = C_\infty$.
	This completes the proof.
\end{pf}
\noindent

\par\smallskip
Take an ultrametric space $(X,d)$ such that the value poset of $d$ is countable and narrow. Then
Theorem~\ref{ThmHevyNarow} applies. Hence $(X,\cB_d)$ is chain intersection stable and
assertions 1) and 2) of Theorem~\ref{MT1} hold. This proves Theorem~\ref{MT2}.

\section{Negative results}               \label{Sectnr}

We will first present an example showing that Theorem~\ref{MT2} is no longer true when the value set is not narrow.
We will need the following simple fact which provides many examples of (generalized) ultrametric spaces.

\begin{prop}                                                 \label{PropuT}
	Take any set $X$ and $\Tau$ a subset of ${\cal P}(X)$ such that for every two distinct elements $x,y\in X$ there is
	a smallest set $B_\Tau(x,y)\in \Tau$ containing $x$ and $y$. Set $u_\Tau(x,y):= B_\Tau(x,y)$ if $x\ne y$ and $u_\Tau(x,x):=
	\emptyset$. Then $u_\Tau$ is an ultrametric on $X$ with value set contained in $\Tau\cup\{\emptyset\}$, which is
	partially ordered by inclusion, with $\bot = \emptyset$. With respect to this ultrametric,
	$$B(x,y) = B_\Tau(x,y).$$
\end{prop}

\begin{pf}
	If $x\ne y$ then by definition, $x,y\in B_\Tau(x,y)$, so $u_\Tau(x,y)\ne\emptyset=\bot$. This together with the
	definition of $u_\Tau(x,x)$ proves (U1). Furthermore, (U3) holds since $B_\Tau(x,y)=B_\Tau(y,x)$.
	In order to prove the ultrametric triangle law (U2), take $x,y,z\in X$ and a set $S\in \Tau$ such that
	$B_\Tau(x,y)\subseteq S$ and $B_\Tau(y,z)\subseteq S$. Then $x,z\in S$ and by the definition of $B_\Tau(x,z)$,
	we find that $B_\Tau(x,z)\subseteq S$.
	
	To prove the second assertion, we observe that
	$$
	z\in B(x,y) \iff u_\Tau(x,z)\leq u_\Tau(x,y) \iff B_\Tau(x,z)\subseteq B_\Tau(x,y).
	$$
	Hence if $z\in B(x,y)$, then $z\in B_\Tau(x,y)$. Conversely, if $z\in B_\Tau(x,y)$, then $x,z\in B_\Tau(x,y)$ and by the
	definition of $B_\Tau(x,z)$, we find that $B_\Tau(x,z)\subseteq B_\Tau(x,y)$, whence $z\in B(x,y)$. This proves that
	$B(x,y) = B_\Tau(x,y)$.
\end{pf}

The following is an immediate consequence of Proposition~\ref{PropuT}.

\begin{wn}
	Each ball space in which for every two distinct elements $x,y$ there is a smallest ball containing $x$ and $y$
	admits a canonical induced ultrametric.
\end{wn}

\begin{exs}
1) If $\cB$ is the ball space of {\it all} closed ultrametric balls in an ultrametric space with linearly ordered
value set, then the value set of the induced ultrametric consists of exactly all precise balls.

\par\smallskip\noindent
2) If $X$ is a $T_1$ topological space, then each two element set is closed. Therefore, if $\cB$ is the ball space
of all closed sets in a $T_1$ topological space, then under the induced ultrametric $u_\Tau\,$, the values $u_\Tau(x,y)$
and $u_\Tau(x',y')$ are incomparable whenever $\{x,y\}\ne \{x',y'\}$. That is, $u_\Tau X^2\setminus\{\bot\}$ is an anti-chain.
If $X$ is not $T_1\,$, then the complexity of the ordering on $u_\Tau X^2 \setminus\{\bot\}$ (given by set inclusion) can
be seen as a measure of how far $X$ is from being $T_1\,$.
\end{exs}

\par\medskip
The next example will prove Theorem~\ref{MT3}:
\begin{ex}                                    \label{Ex1}
Take $X = (\omega+1) \times \omega$. For $m,k < \omega$, we set
$$B_{m,k} := \setof{(n,k)}{m \loe n \loe \omega} \cup \setof{(\omega, \ell)}{k \loe \ell < \omega}$$
and
$$\Tau := \setof{B_{m,k}}{m,k \in \omega} \cup \setof{\dn{(\omega, k)}{(\omega, \ell)}}{k, \ell \in \omega, \; k \ne \ell} \cup \sn X.$$
Then for all distinct $(m, k), (n, \ell) \in X$ with $m \loe n$, there is a smallest set in $\Tau$ containing them:
	\begin{equation}
	B_\Tau((m,k),(n,\ell)) = \left\{
	\begin{array}{ll}
		B_{m,k}& \mbox{if } m\leq n<\omega \mbox{ and } k=\ell\,, \\
		& \mbox{or } m<n=\omega \mbox{ and } k\leq\ell\\[.2cm]
		\{(\omega,k),(\omega,\ell)\} & \mbox{if } m=n=\omega \mbox{ and } k\ne\ell\\[.2cm]
		X & \mbox{if } m\leq n<\omega \mbox{ and } k\ne\ell\,, \\
		& \mbox{or } m<n=\omega \mbox{ and } k>\ell\,, \\
		& \mbox{or }  m=n=\omega \mbox{ and } k=\ell\>.
	\end{array}
	\right.
	\end{equation}
In view of the symmetry of the sets $B_\Tau$, the assumption $m \loe n$ is no loss of generality.
	
By Proposition~\ref{PropuT} we obtain an induced ultrametric $u_\Tau$, and the family $\Bee_\Tau$
of its closed ultrametric balls is equal to $\Tau$.
We show that $(X, u_\Tau)$ is spherically complete.
	
Take a chain $\calN$ in $\Tau$ . If
$\calN$ contains a smallest ball, then its intersection is equal to this ball and
hence nonempty. If $\calN$ does not contain a smallest ball, then it must be of the form
$\setof{B_{m_i,k}}{i \in \omega}$ for some $k \in \omega$, where $(m_i)_{i \in \omega}$ is a strictly increasing sequence in $\omega$. The intersection of this chain is the nonempty set $\setof{(\omega, \ell)}{\ell \goe k}$.
	
Suppose that $(X, \Bee')$ is a chain intersection closed expansion of $(X, \Bee_\Tau)$. Then by what we have just 
shown,
$$\calN = \setof{ \setof{(\omega,\ell)}{k \loe \ell < \omega}}{k \in \omega} \subs \Bee'.$$
We find that $(X, \Bee')$ is not spherically complete, since $\calN$ is a nest and its
intersection is empty.

Let us note that the arguments of the proof show that with $\cB'=\cB_\Tau\cup\calN$,
the expansion $(X,\cB')$ of $(X,\cB_\Tau)$ is chain intersection closed. It follows that $\cB'=\ci(\cB)$, showing that while assertion 2) of Theorem~\ref{MT1} fails, assertion 1) still holds.
\end{ex}

Our final result shows that the assertion of Theorem~\ref{ThmHevyNarow} does not remain true when the 
countability assumption is dropped.
\begin{tw}                                     \label{MT4}
There exists an uncountable spherically complete narrow ball space, closed under finite intersections, which is not 
chain intersection stable and does not admit any expansion that is chain intersection closed and spherically 
complete.
\end{tw}
\noindent
This theorem is proved by the following example.
\begin{ex}                                    \label{Ex2}
Let $X = ((\omega_1 + 1) \times (\omega + 1)) \setminus \sn{(\omega_1,\omega)}$, so that the elements of $X$ are pairs of ordinals $(\al,\beta)$, where $\al \loe \omega_1$, $\beta \loe \omega$ and $(\al,\beta) \ne (\omega_1, \omega)$.
	Define $B_{\al, n} = \setof{(x,y) \in X}{\al \loe x, \; n \loe y}$. This is the maximal rectangle in $X$, whose bottom-left vertex is $(\al,n)$.
	Let $$\Bee = \setof{B_{\al,n}}{\al<\omega_1, \; n < \omega}.$$
	Then $(X, \Bee)$ is a ball space which is closed under finite intersections, because $B_{\al,n} \cap B_{\beta, m} = B_{\gamma, k}$ where $\gamma = \max(\al,\beta)$ and $k = \max(n,m)$.
	It is clear that $(\Bee, \subs)$ is narrow, as it is isomorphic to $\omega_1 \times \omega$ with the product ordering.

If $\Cee$ is a chain in $\Bee$ then $\bigcap \Cee \nnempty$. Indeed, this is true if $\Cee$ is countable because then $[\al,\omega) \subs \bigcap \Cee$ for a sufficiently big $\al <\omega_1$. Otherwise, if for each $\al < \omega_1$ there is $\xi_\al > \al$ such that $B_{\xi_\al,n_\al} \in \Cee$ for some $n_\al < \omega$, then there are uncountably many $\al$s such that $n_\al = n$ is constant and consequently $(\omega_1, n) \in \bigcap \Cee$, because $\Cee$ is a chain. This proves that $\Bee$ is spherically complete.

On the other hand, $(X, \Bee)$ is not chain intersection stable, because it contains the chains $\Cee_\al = 
\sett{B_{\al,n}}{n < \omega}$ for which $\bigcap \Cee_\al = ([\al,\omega_1) \times \sn \omega) \cap X$ and hence 
$\bigcap_{\al < \omega_1} \bigcap \Cee_\al = \emptyset$. Every chain intersection closed expansion of $\Bee$ will 
contain the chain $\sett{\bigcap \Cee_\al}{\al < \omega_1}$ and will thus not be spherically complete.	
\end{ex}

Note that in contrast to Example~\ref{Ex1}, Proposition~\ref{PropuT} cannot be applied here to derive an 
ultrametric space, because for any two points whose second coordinate is $\omega$ there is no smallest ball 
in $\Bee$ containing them.

\end{document}